\documentclass{aptpub}

\authornames{Harry Crane and Peter McCullagh} 
\shorttitle{Reversible Markov structures on divisible set partitions} 

\def\ESF{\mathop{\rm ESF}\nolimits}

\def\crp{\mathop{\rm CRP}\nolimits}

\def\PY{\mathop \text{Ewens--Pitman}\nolimits}

\def\integerpartitions{\mathop{\mathcal{P}_n}\nolimits}

\def\evenpartitionsnk{\mathop{\mathcal{P}_{[nk]\,\mid\, k}}\nolimits}
\def\evenpartitionsnkm{\mathop{\mathcal{P}_{[nk]\,\mid\, k}^{(m)}}\nolimits}

\def\partitionsNm{\mathop{\mathcal{P}_{\mathbb{N}}^{(m)}}\nolimits}
\def\partitionsn{\mathop{\mathcal{P}_{[n]}}\nolimits}
\def\partitionsnm{\mathop{\mathcal{P}_{[n]}^{(m)}}\nolimits}

\def\partitionsN{\mathop{\mathcal{P}_{\mathbb{N}}}\nolimits}

\def\part{\mathop{\mbox{Part}}\nolimits}

\def\ell{\mathop{[l]}\nolimits}

\begin{document}

\title{Reversible Markov structures on divisible set partitions} 

\authorone[Rutgers University]{Harry Crane} 

\addressone{Rutgers University\\ Department
of Statistics \\ 
110 Frelinghuysen Road\\
Piscataway, NJ 08854, USA} 

\authortwo[University of Chicago]{Peter McCullagh}

\addresstwo{University of Chicago\\
Department of Statistics\\
Eckhart Hall\\
5734 S.\ University Ave.\\
Chicago, IL 60637}

\begin{abstract}
We study $k$-divisible partition structures, which are families of random set partitions whose block sizes are divisible by an integer $k=1,2,\ldots$. 
In this setting, exchangeability corresponds to the usual invariance under relabeling by arbitrary permutations; however,
for $k>1$, the ordinary deletion maps on partitions no longer preserve divisibility, and so a random deletion procedure is needed to obtain a partition structure.  
We describe explicit Chinese restaurant-type seating rules for generating families of exchangeable $k$-divisible partitions that are consistent under random deletion.
We further introduce the notion of {\em Markovian partition structures}, which are ensembles of exchangeable Markov chains
on $k$-divisible partitions that are consistent under a random process of {\em Markovian deletion}.
The Markov chains we study are reversible and refine the class of Markov chains introduced in {\em J.\ Appl.\ Probab.}~{\bf48}(3):778--791.
\end{abstract}



\section{Introduction}
A {\em partition} $\pi$ of $[n]:=\{1,\ldots,n\}$ is a collection $B_1/\cdots/B_m$
of non-empty, disjoint subsets, called blocks, for which $\bigcup_{1\leq j\leq m}B_j=[n]$.
For $n,k\geq1$, we call a partition $\pi=B_1/\cdots/B_m$ of $[nk]$
{\em $k$-divisible}, or just {\em divisible}, if the cardinality of each block $B_1,\ldots,B_m$
is divisible by $k$.  When $k=2$, we call $\pi$ an {\em even partition}.  

Divisible partitions are natural in ecological applications as well as 
randomization in experimental design.  For example, in experimental design,
each of $nk$ individuals is assigned one of $k$ treatments.  If individuals are
further grouped into blocks so that every treatment is assigned the same number of
times within each block, then the block structure of the design is a divisible partition of $[nk]$.  In this
setting, divisible partitions are related to group-divisible association schemes; see Bailey \cite{BaileyAssociationSchemes}
for further connections between experimental design and the theory of partitions.
Our study of probabilistic structures of divisible random partitions is motivated by the above
heuristic as well as the appeal of partition models to applications in clustering and classification \cite{Crane2014cluster,McCullaghYang2008},
population genetics \cite{Ewens1972, FisherCorbetWilliams1943}, and linguistics \cite{EfronThisted1976,EfronThisted1987}.

We study ensembles of random divisible partitions whose distributions are consistent under a random deletion operation.   
Kingman \cite{Kingman1978a, Kingman1978b} first studied the deletion properties of random integer partitions.  He defined a {\em
partition structure} as a collection $\mathbb{P}:=(\mathbb{P}_1,\mathbb{P}_2,\ldots)$ of probability distributions on the spaces $(\integerpartitions,\,n\in\mathbb{N})$ of
integer partitions of each $n=1,2,\ldots$ such that if
$\lambda'\in\mathcal{P}_{n-1}$ is obtained by choosing a part of
$\lambda\sim\mathbb{P}_n$ with probability proportional to its size and
reducing it by one, then $\lambda'\sim\mathbb{P}_{n-1}$.  On set
partitions, consistency of a family of distributions $\mathbb{P}$ is defined through the usual non-random restriction operation.
Any consistent collection determines a unique probability measure on partitions of $\mathbb{N}:=\{1,2,\ldots\}$.  Gnedin, et al  \cite{GnedinHaulkPitman} study further
deletion properties of random partitions. 

In our main theorems, we extend the Chinese restaurant process \cite{Pitman2005} and the Ewens--Pitman Markov chain \cite{Crane2011a} to the space of $k$-divisible partitions, and we obtain the finite-dimensional distributions of these processes.  By reversing these procedures, we describe natural deletion mechanisms under which the prescribed finite-dimensional distributions are consistent.
In the Markov chain case, this produces a {\em Markovian partition structure},
that is, a family $\{(\varepsilon^{(n)},\mathcal{E}^{(n)})\}_{n\geq1}$ of initial distributions $\varepsilon^{(n)}$
and transition probability measures $\mathcal{E}^{(n)}$ on $k$-divisible partitions of $[nk]$ that are
exchangeable and consistent under a Markovian deletion scheme. 
Under this operation, Markovian partition structures have the form
\begin{eqnarray*}
\mathcal{P}_{[nk]\,\mid\, k}&\quad\longrightarrow\quad&\mathcal{P}_{[nk]\,\mid\, k}\\
\bigg\downarrow & &\bigg\downarrow\\
\mathcal{P}_{[(n-1)k]\,\mid\, k} & \quad\longrightarrow\quad & \mathcal{P}_{[(n-1)k]\,\mid\, k},
\end{eqnarray*}
where horizontal arrows denote Markov transitions in time and vertical arrows represent randomized projections by the Markovian deletion scheme.
In particular, the marginal distributions at fixed times, as $n$ varies, are consistent; and the marginal distribution of each sequence, for fixed $n$, is a Markov chain.
For a single time $t$, the marginal behavior of the Markovian deletion scheme coincides with the deletion scheme for $k$-divisible partition structures.
We provide the details in Sections \ref{section:even partition structures} and \ref{section:Markovian deletion}.

To obtain  a characteristic measure
$\varepsilon^{(\infty)}$ on the limit space of partitions of $\mathbb{N}$, we need a  deterministic deletion operation; but there is no such operation in the $k$-divisible setting with $k>1$.
When $k>1$, simple deletion of the highest labeled group $\{nk+1,\ldots,(n+1)k\}$ from an exchangeable divisible partition
of $[(n+1)k]$ does not preserve divisibility; therefore, a random deletion scheme is needed.  
For example, the partition $\pi=1468/27/35$ is $2$-divisible, but the partition $\pi'=146/2/35$ obtained by deleting elements $\{7,8\}$ is not $2$-divisible because the cardinalities of $\{1,4,6\}$ and $\{2\}$ are not even.

Though the processes we study are not sampling consistent in the ordinary sense, the finite-dimensional
processes we generate are exchangeable and, in the Markov chain case, reversible with respect to the $k$-divisible extension of the Ewens--Pitman distribution \cite{Ewens1972, Pitman2005}.
Reversible processes for partition-valued Markov chains have been studied previously; see \cite{Bertoin2008,Crane2011a}.

Our main discussion focuses on the analog to the Ewens distribution for divisible partitions and Markov
chains; however, these conclusions apply more generally to any paintbox measure.
 We develop these ideas formally in Sections
\ref{section:even partition structures} and \ref{section:even Markov structures}. 
We make some concluding remarks in Section \ref{section:remarks}.

\section{Preliminaries: Random partitions}
A partition $\mathbf{n}$ of  $n=1,2,\ldots$ is a non-increasing collection of positive integers $\mathbf n = (n_1, \ldots, n_k)$, called {\em parts} or {\em summands}, whose sum is~$n$.
Alternatively, $\mathbf{n}$ can be expressed in terms of its multiplicities
$(\lambda_1,\ldots, \lambda_n)$,
$\mathbf n = 1^{\lambda_1} 2^{\lambda_2}\cdots n^{\lambda_n}$,
such that $\sum\lambda_j = k$ is the number of parts
and $\sum j \lambda_j = n$.
%
For each
$n\in\mathbb{N}$, the Ewens sampling formula  with parameter
$\theta >0$, ESF($\theta$), is the probability distribution on integer partitions of $n$ with closed form expression
\begin{equation}\label{eq:ESF}
p_n(\mathbf{n};\theta)=\frac{\theta^{\#\mathbf{n}}}{\theta^{\uparrow
n}}\frac{n!}{\prod_{j=1}^n
j^{\lambda_j}\lambda_j!},\quad\mathbf{n}=1^{\lambda_1}\cdots n^{\lambda_n}\in\integerpartitions,
\end{equation}
where $\#\mathbf{n}:=\sum_{j=1}^n\lambda_j$ is the number of parts of $\mathbf{n}$
and $\theta^{\uparrow n}:=\theta(\theta+1)\cdots(\theta+n-1)$.  
Ewens \cite{Ewens1972} first derived \eqref{eq:ESF} while studying the sampling theory of neutral alleles, but the Ewens sampling formula also occurs in purely mathematical contexts, e.g., as the asymptotic distribution of large prime factors \cite{DonnellyGrimmett1993} and as a special case of the $\alpha$-permanent of a matrix \cite{Crane2013b,Crane2013LAA}.

Ewens's sampling formula more naturally resides on the space
$\partitionsn$ of partitions of $[n]$, where it is the $(0,\theta)$ sub-family of the two-parameter $\PY(\alpha,\theta)$ family with
finite-dimensional marginal distributions
\begin{equation}\label{eq:fidi PY}
\varrho_{\alpha,\theta}^{(n)}(\pi)=\frac{(\theta/\alpha)^{\uparrow\#\pi}}{\theta^{\uparrow
n}}\prod_{b\in\pi}-(-\alpha)^{\uparrow\#b},\quad\pi\in\partitionsn,\end{equation}
where $\#\pi$ denotes the number of blocks of $\pi$, $\#b$ the cardinality of block $b$, and
$(\alpha,\theta)$ satisfies either
\begin{itemize}
        \item $\alpha=-\kappa<0$ and $\theta=m\kappa$, for $m=1,2\ldots$, or
        \item $0\leq\alpha\leq1$ and $\theta>-\alpha$.
\end{itemize}
 For fixed $(\alpha,\theta)$, $\varrho_{\alpha,\theta}:=(\varrho_{\alpha,\theta}^{(n)},\,n\in\mathbb{N})$
is a consistent collection of exchangeable probability measures on the system $(\partitionsn,\,n\in\mathbb{N})$.  
In particular, $(\varrho^{(n)}_{\alpha,\theta},\,n\in\mathbb{N})$ is
\begin{itemize}
	\item {\em exchangeable} if, for every $n\in\mathbb{N}$, $\varrho^{(n)}_{\alpha,\theta}(\pi)$ depends on $\pi$ only through its block sizes, and
	\item {\em consistent under subsampling} if, for all $m\leq n$, the image measure of $\varrho_{\alpha,\theta}^{(n)}$ by restriction to $\mathcal{P}_{[m]}$ is $\varrho_{\alpha,\theta}^{(m)}$; that is, $\varrho_{\alpha,\theta}^{(m)}=\varrho_{\alpha,\theta}^{(n)}\mathbf{R}_{m,n}^{-1}$, where $\mathbf{R}_{m,n}:\partitionsn\rightarrow\mathcal{P}_{[m]}$ denotes the restriction map,
\[\mathbf{R}_{m,n}\pi:=\{B_1\cap[m],\ldots,B_l\cap[m]\}\setminus\{\emptyset\},\quad\pi=B_1/\cdots/B_l\in\partitionsn.\]
\end{itemize}

As a result, the finite-dimensional marginals in \eqref{eq:fidi PY} determine a
unique probability measure $\varrho_{\alpha,\theta}$ on $\partitionsN$, called the $\PY(\alpha,\theta)$ distribution.  Throughout the paper, we assume that a pair $(\alpha,\theta)$ is always within the parameter space of the
Ewens--Pitman model, and we call any random partition with finite-dimensional distributions \eqref{eq:fidi PY} an {\em $(\alpha,\theta)$-partition}. 
The distribution on $\partitionsn$ corresponding to $\ESF(\theta)$ in
\eqref{eq:ESF} is $\PY(0,\theta)$.

Sampling consistency of the Ewens--Pitman family is easily observed through its  {\em Chinese restaurant} construction.  
We construct a sequence $\Pi:=(\Pi_1,\Pi_2,\ldots)$ of finite exchangeable partitions by putting
$\Pi_1=\{1\}$ and, given $\Pi_n=\pi=B_1/\cdots/B_m\in\partitionsn$, we generate $\Pi_{n+1}$ by inserting the element $n+1$ 
\begin{itemize}
        \item in occupied block $b\in\pi$ with probability $(\#b-\alpha)/(n+\theta)$ and
        \item in its own block of $\Pi_{n+1}$ with probability $(\theta+m\alpha)/(n+\theta)$.
\end{itemize}
We write $\crp(n,\alpha,\theta)$ to denote the conditional probability distribution
of $\Pi_{n+1}$ given $\Pi_n$ above.

More generally, Kingman's paintbox process \cite{Kingman1978a} describes the law of any exchangeable partition of $\mathbb{N}$.  Let $\nu$ be a probability measure
on
\[\Delta^{\downarrow}:=\left\{(s_1,s_2,\ldots):\,s_1\geq s_2\geq \cdots\geq0,\,\sum_{j}s_j\leq 1\right\}.\]
Given $(S_1,S_2,\ldots)\sim\nu$ and $S_0:=1-\sum_{j}S_j$, we generate a sequence $X:=(X_1,X_2,\ldots)$ of conditionally independent
random variables with
\[P\{X_i=j\mid S\}=\left\{\begin{array}{cc}
S_j,& j\geq1\\
S_0,& j=-i\\
0,& \text{otherwise.}
\end{array}\right.\]
The partition $\Pi$ of $\mathbb{N}$ defined by
\[i\text{ and }j\text{ are in the same block of }\Pi\quad\text{if and only if}\quad X_i=X_j\]
is exchangeable and obeys {\em (Kingman's) paintbox distribution}, or {\em paintbox process}, directed by $\nu$, denoted $\Pi\sim\varrho_{\nu}$.  The $\PY(\alpha,\theta)$ law
corresponds to the paintbox process directed by the $\text{Poisson--Dirichlet}(\alpha,\theta)$ law on $\Delta^{\downarrow}$.

In the next section, we give Chinese restaurant-type constructions for divisible
partitions.   Through this process, we construct a partition structure
on exchangeable $k$-divisible partitions of $[nk]$.  We also introduce Markovian partition structures, which are families of Markov chains consistent under a random {\em Markovian} deletion process.  This extends Kingman's partition structures to sequences of exchangeable $k$-divisible partitions and further refines the family of exchangeable Markov chains studied in \cite{Crane2011a}.

\section{Divisible partitions}\label{section:even partitions}
We call any subset $A\subseteq\mathbb{N}$ {\em $k$-divisible}, or just {\em divisible}, if $\#A$ is divisible by $k$, denoted $A\mid k$.
We call a partition $\pi$ {\em $k$-divisible} if each of its blocks is $k$-divisible, denoted $\pi\mid k$, and we write $\evenpartitionsnk$ to denote the space of
$k$-divisible partitions of $[nk]$.  For example, $\pi=148/235679$ is $k$-divisible
for $k=1$ and $k=3$.

\subsection{Chinese restaurant construction for divisible partitions}\label{section:CRP
even}
Let $k\in\mathbb{N}$, $(\alpha,\theta)$ be fixed, and $\mathbb{N}$ be a
population of individuals.  Suppose individuals arrive at a restaurant in
groups of size $k$, regarded as $\{ik+1,\ldots,ik+k\}$ for $i=1,2,\ldots$.  In
particular, for every $n\in\mathbb{N}$, we construct a random divisible partition
$\Pi_n\in\evenpartitionsnk$ according to the following seating rule.  As usual, the
tables in the restaurant correspond to the blocks of a random partition.

\vspace{3mm}
\noindent{\bf Divisible random seating rule} 
\begin{itemize}
        \item[(1)] The first $k$ individuals are seated at the same table,
$\Pi_1:=\{\{1,\ldots,k\}\}$.
        \item[(2)] After $nk$ individuals are seated according to $\Pi_n$, the next $k$
individuals $nk+1,\ldots,nk+k$ are seated randomly as follows.  We initialize by putting $\Pi_n^{(1)}=\Pi_n$.
        \begin{itemize}
                \item[(a)] Independently for each $i=2,\ldots,k$, $nk+i$ chooses $u^{(i)}$
uniformly among $[nk+i-1]$ and immediately displaces $u^{(i)}$ in $\Pi_n^{(i-1)}$ to define $\Pi_n^{(i)}$.  If the chosen element $u^{(i)}$ is not in $\Pi_n^{(i-1)}$, then no displacement occurs.
                \item[(b)] After each individual has made its choice of $u^{(i)}$ above, there
are $k$ individuals $(nk+1,w^{(2)},\ldots,w^{(k)})$ waiting to be seated. 
The group $w^*:=\{nk+1,w^{(2)},\ldots,w^{(k)}\}$ is treated as a single unit and randomly
chooses a table $b\in\Pi_n^{(k)}\cup\{\emptyset\}$ according to
$\crp(nk,\alpha,\theta)$, that is,
\[P\{w^*\mapsto b\mid \Pi_n^{(k)}=\pi\}\propto\left\{\begin{array}{cc} \#b-\alpha,& b\in\pi\\ \alpha\#\pi+\theta,& b=\emptyset.\end{array}\right.\]
We define $\Pi_{n+1}$ as the partition resulting from this seating assignment.
        \end{itemize}

        \end{itemize}

This construction generates a collection $(\Pi_1,\Pi_2,\ldots)$ of divisible partitions of $[nk]$. 

\begin{remark}
Note the change in notation in (b) above from $u^{(i)}$, denoting displaced individuals,
to $w^{(i)}$, denoting those individuals still not seated after the final displacement.
This reflects the possibility that a single element can be displaced multiple times, because a particular element can be chosen as $u^{(i)}$ for multiple $i=2,\ldots,k$.
The following example illustrates the seating procedure.
\end{remark}

\begin{ex}\label{ex:CRP}
With $n=k=3$ and fixed $(\alpha,\theta)$, we generate the partition $\Pi_3=134689/257$ from the above seating procedure as follows.
\begin{itemize}
	\item We begin with $\Pi_1=123$;
	\item At time 2, individuals $4,5,6$ arrive and, following Step (2a), element 5 first chooses $u^{(2)}$ from $\{1,2,3,4\}$, say $u^{(2)}=4$, and element 6 chooses $u^{(3)}$ from $\{1,2,3,4,5\}$, say $u^{(3)}=2$.  Then, after Step 2, we have partition $136$, with $2,4,5$ displaced.  (Note that the displaced individuals differ from the set $\{4,u^{(2)},u^{(3)}\}=\{2,4\}$.)  Treating $w^*:=\{2,4,5\}$ as a single unit, we choose to put $\{2,4,5\}$ in the same block as $\{1,3,6\}$ with probability $(3-\alpha)/(3+\theta)$ to obtain $\Pi_2=123456$.
	\item At time 3, individuals $7,8,9$ arrive and choose $u^{(2)}=2$ and $u^{(3)}=5$ so that the partition after Step (2a) is $134689$.  We now place $w^*=\{2,5,7\}$ in its own block with probability $(\alpha+\theta)/(6+\theta)$ to get $\Pi_3=134689/257$.
\end{itemize}
Note that there is more than one way to obtain $\Pi_3=134689/257$.  
We derive the distribution of $\Pi_n$ in Theorem \ref{thm:CRP even} below.

\end{ex}

In proving the following theorem, and throughout the paper, we write $i\sim_{\pi}j$ to denote that elements $i$ and $j$ are in the same block of a partition $\pi$.
Exchangeability on $\evenpartitionsnk$ is defined in the usual way: $\varepsilon^{(n)}_{\alpha,\theta}$ is exchangeable if
$\varepsilon^{(n)}_{\alpha,\theta}(\pi)$ depends only on $(\#b,\,b\in\pi)$, for every $\pi\in\evenpartitionsnk$.

        \begin{theorem}\label{thm:CRP even}
        Let $\Pi:=(\Pi_1,\Pi_2,\ldots)$ be a sequence of random partitions generated by the
above seating rule.  Then each $\Pi_n$ is marginally an exchangeable $k$-divisible partition of $[nk]$ with
distribution 
 \begin{equation}\label{eq:even fidi}
\varepsilon^{(n)}_{\alpha,\theta}(\pi)=\frac{n!}{(nk)!}\frac{(\theta/\alpha)^{\uparrow\#\pi}}{(\theta/k)^{\uparrow
n}} \prod_{b\in\pi}-(-\alpha/k)^{\uparrow(\#b/k)}\frac{\#b!}{(\#b/k)!},
\quad\pi\in\evenpartitionsnk.\end{equation}
        
\end{theorem}
\begin{proof}
To establish \eqref{eq:even fidi}, we fix $k\in\mathbb{N}$ and $(\alpha,\theta)$ in the parameter space of the Ewens--Pitman model.  We show \eqref{eq:even fidi} by
induction on $n$. Clearly, \eqref{eq:even fidi} holds for $n=1$ since
$\varepsilon^{(1)}_{\alpha,\theta}$ in \eqref{eq:even fidi} is the point mass at the
single-block partition $\{1,2,\ldots,k\}$.

Now, assume that \eqref{eq:even fidi} holds for $n\in\mathbb{N}$ and consider
$\pi^*\in\mathcal{P}_{[(n+1)k]\,\mid\, k}$.  Let
$A_{\pi^*}:=\{\pi\in\evenpartitionsnk:\pi\mapsto\pi^*\}$ be divisible partitions
of $[nk]$ for which there is a positive probability of the event $\{(\Pi_n,\Pi_{n+1})=(\pi,\pi^*)\}$ in the divisible seating process.  By definition, the block sizes of each
$\pi\in A_{\pi^*}$ and $\pi^*$ are identical except for the block to which the displaced
group $(nk+1,w^{(2)},\ldots,w^{(k)})$ is inserted in $\pi$ during step (2b).  Let this block be written as $b^*\in\pi^*$.

Each random displacement in step (2a) has probability $1/[(nk+1)\cdots(nk+k-1)]$ and the random table assignment in step (2b)
follows $\crp(nk,\alpha,\theta)$, which assigns probability
$(\#b-\alpha)/(nk+\theta)$ to $b$ if $b\neq\emptyset$ and
$(\theta+\alpha\#\pi)/(nk+\theta)$ if $b=\emptyset.$  By the induction hypothesis, each $\pi\in A_{\pi^*}$ is distributed as in \eqref{eq:even fidi}, so that the joint probability of the event $(\Pi_n,\Pi_{n+1})=(\pi,\pi^*)$, for every $\pi\in A_{\pi^*}$, is
\[\frac{n!}{(nk)!}\frac{(\theta/\alpha)^{\uparrow\#\pi}}{(\theta/k)^{\uparrow
n}} \prod_{b^*\in\pi^*}-(-\alpha/k)^{\uparrow(\#b^*/k)}\prod_{b\in\pi}\frac{\#b!}{(\#b/k)!},\]
for every $n\geq1$.

Since this joint probability is the same for all pairs $(\pi,\pi^*)$, we obtain the marginal distribution of $\Pi_{n+1}$ by multiplying the number of partitions in $A_{\pi^*}$.
Let $\overleftarrow{\pi}_n^*$ denote the labeled $(nk+1)$-shift of $\pi^*$.  That
is, define $\overleftarrow{\pi}^*_n$ as a partition of $\{2,\ldots,n\}$ with
\[i\sim_{\overleftarrow{\pi}^*_n} j\quad\text{if and only if}\quad nk+i\sim_{\pi^*} nk+j,\]
and label each $b\in\overleftarrow{\pi}^*_n$ by its smallest element.
In step (2a) of the random seating plan, $\overleftarrow{\pi}^*_n$ is obtained by
random displacement of the elements of $\pi$.  Let $\sigma^*$ be the permuation $[nk]\rightarrow[nk]$ defined by the product of transpositions
\[\sigma^*:=\begin{pmatrix} nk+2 \\ u^{(2)}\end{pmatrix}\cdots\begin{pmatrix} nk+k \\ u^{(k)}\end{pmatrix}\]
and let $\varphi_{\pi,\pi^*}:\{2,\ldots,n\}\rightarrow[(n+1)k]$ be the operation
corresponding to step (2a) of the above divisible seating rule: for $i=2,\ldots,k$,
$\varphi_{\pi,\pi'}(i)=\sigma^{*}(nk+i)$ is the element occupying position $i$ in $(nk+1,w^{(2)},\ldots,w^{(k)})$.  Let $\varphi_{\pi,\pi^*}(\overleftarrow{\pi}^*_n)$
denote the partition obtained by replacing each $i\in\{2,\ldots,n\}$ with
$\varphi_{\pi,\pi^*}(i)$.  Writing $b\in\pi$ to denote the block to which the
displaced group is added in step (2b) of the seating process, we have
\[{{\#b+k-1}\choose{k-1}}=\frac{(\#b+k-1)!}{(k-1)!\#b!}\]
choices of the elements of $\varphi_{\pi,\pi^*}(\overleftarrow{\pi}^*_n)$ for
any choice $(u^{(2)},\ldots,u^{(k)})$ of displaced elements; and there are
\[\frac{(k-1)!}{\prod_{b'\in\overleftarrow{\pi}^*_n}\#b'!}\]
ways to arrange these elements into a labeled partition with block sizes
corresponding to the block sizes of $\overleftarrow{\pi}^*_n$.  Finally, the
assignments under $\varphi_{\pi,\pi^*}$ within each $b'\in\overleftarrow{\pi}_n^*$ can be
rearranged in $\#b'!$ ways to obtain a total of $(\#b+k-1)\cdots(\#b+1)$ partitions in $A_{\pi^*}$.
Equation \eqref{eq:even fidi} now follows by the induction hypothesis.
This completes the proof.

\end{proof}

\begin{remark}
Note the relationship between distributions \eqref{eq:fidi PY} and \eqref{eq:even fidi}.
These distributions coincide when $k=1$, because (2a) of the divisible seating rule is nugatory and the divisible random seating rule is equivalent to the usual Chinese restaurant seating rule in this case.
Otherwise, these distributions differ as a result of the random shuffling that occurs during step (2a) of the divisible random seating rule.

\end{remark}

For $\alpha=0$, \eqref{eq:even fidi} with parameter $(0,\theta)$ becomes 
\begin{equation}\label{eq:ewens even fidi}
\varepsilon^{(n)}_{0,\theta}(\pi)=\frac{n!}{(nk)!}\frac{\theta^{\#\pi}\prod_{b\in\pi}(\#b-1)!}{(\theta/k)^{\uparrow
n}},\end{equation}
the marginal distribution of each $\Pi_n$ in the construction with random seating
rule $\crp(nk,0,\theta)$ in step (2b).  From \eqref{eq:ewens even fidi}, we obtain
the combinatorial identity
\begin{equation}\label{eq:identity}\frac{1}{(nk)!}\sum_{\pi\in\evenpartitionsnk}\theta^{\#\pi}\Gamma(\pi)=\left.\left(\frac{\theta}{k}\right)^{\uparrow
n}\right/ n!,\end{equation}
where $\Gamma(\pi):=\prod_{b\in\pi}(\#b-1)!$. 
This identity gives the generating function for {\em $k$-divisible permutations}, 
that is, permutations of $[nk]$ whose cycle sizes are all divisible by $k$, and also gives the following special property of the Ewens distribution.

\begin{cor}\label{cor:conditioned ewens even}
For $\theta>0$, the distribution $\varepsilon^{(n)}_{0,\theta}$ in \eqref{eq:even
fidi} is that of a $(0,\theta)$-partition conditioned to be $k$-divisible.
\end{cor}

\subsection{Divisible partition structures}\label{section:even partition structures}
We specify a random deletion scheme for
$\evenpartitionsnk$ as follows.  Given $\pi\in\evenpartitionsnk$, let $b_*\in\pi$
denote the block of $\pi$ containing $(n-1)k+1$.  
\begin{itemize}
	\item[(i)] Sequentially, for $i=k,k-1,\ldots,2$, an element $u^{(i)}$ is chosen uniformly from the set
	\[(b_*\cap[(n-1)k+i])\setminus\{(n-1)k+1,u^{(k)},\ldots,u^{(i+1)}\}.\]
 Let $\pi^*$ be the image of $\pi$ under permutation by the product of transpositions
\begin{equation}\label{eq:even displacement}
\sigma^{*}:=\begin{pmatrix} (n-1)k+2 \\ u^{(2)}\end{pmatrix} \cdots \begin{pmatrix}(n-1)k+k \\ u^{(k)}\end{pmatrix}.\end{equation}
  \item[(ii)] Obtain $\pi'\in\mathcal{P}_{[(n-1)k]\,\mid\, k}$ by deleting
$\{(n-1)k+1,\ldots,(n-1)k+k\}$ from $\pi^*$.  
\end{itemize}

\begin{definition}
We call a collection $\varepsilon=(\varepsilon^{(n)},\,n\in\mathbb{N})$ of probability distributions a {\em divisible partition structure} if, for every $n\in\mathbb{N}$, $\varepsilon^{(n)}$ is the distribution of $\Pi'$ obtained by applying the divisible deletion scheme to $\Pi\sim\varepsilon^{(n+1)}$.
\end{definition}

\begin{theorem}\label{prop:even partition structures}
For any $(\alpha,\theta)$, $\varepsilon_{\alpha,\theta}:=(\varepsilon^{(n)}_{\alpha,\theta},\,n\in\mathbb{N})$ is
a divisible partition structure under the above deletion scheme.\end{theorem}
\begin{proof}
For $\pi'\in\mathcal{P}_{[(n-1)k]\,\mid\, k}$, we define
$A_{\pi'}:=\{\pi\in\evenpartitionsnk:\pi'\leftarrow\pi\}$ to be the set of divisible
partitions of $[nk]$ for which there is positive probability of obtaining $\pi'$
from the above deletion scheme.  Any $\pi\in A_{\pi'}$ has the same block structure
as $\pi'$ except for the block $b_*\in\pi$ containing $(n-1)k+1$, which has $k$ more
elements than its corresponding block in $\pi'$ and will be reduced by $k$ during
the deletion process.  From the proof of Theorem \ref{thm:CRP even}, we can
express the probability of $\pi\in A_{\pi'}$ as
\begin{equation}\label{eq:expression1}
\varepsilon_{\alpha,\theta}^{(n)}(\pi)=\varepsilon_{\alpha,\theta}^{(n-1)}(\pi')\frac{(\#b_*-1)^{\downarrow(k-1)}}{(nk-1)^{\downarrow(k-1)}}\left[\frac{\#b_*-k-\alpha}{(n-1)k+\theta}\mathbb{I}_{\{\#b_*>k\}}+\frac{\theta+\alpha\#\pi'}{(n-1)k+\theta}\mathbb{I}_{\{\#b_*=k\}}\right].
\end{equation}
For any $b\in\pi'$, the probability that $(n-1)k+i$, $i=2,\ldots,k$, displaces a
specific element of $b_*$ is $1/(b_*-k+i-1)$ and, in total, there are $nk-k+i-1$ elements
which $(n-1)k+i$ has the option of displacing.  (Every $\pi\in A_{\pi'}$ corresponds to a choice $b\in\pi'$ to insert the displaced group in step (2b) of the random seating rule.  Given $b\in\pi$, the choice $((n-1)k+1,u^{(2)},\ldots,u^{(k)})$ corresponds to a unique $k$-tuple of transpositions $\sigma^*$ in \eqref{eq:even displacement} to obtain $\pi'$ from $\pi$ by the deletion process.  There are $(nk-k+1)\cdots(nk+1)$ total choices for every such $b\in\pi$.)  By the law of cases, we have
\begin{eqnarray*}
\mathbb{P}\{\Pi_{n-1}=\pi'\}&=&\sum_{\pi\in A_{\pi'}}\mathbb{P}\{\Pi_{n-1}=\pi'\,|\,\Pi_n=\pi\}\varepsilon^{(n)}_{\alpha,\theta}(\pi)\\
&=&\sum_{b_*}\sum_{\sigma^*}\varepsilon^{(n)}_{\alpha,\theta}(\pi)\frac{1}{(\#b_*-1)\cdots(\#b_*-k+1)}\\
&=&\varepsilon^{(n-1)}_{\alpha,\theta}(\pi')\sum_{b_*}\left[\frac{\#b_*-k-\alpha}{(n-1)k+\theta}\mathbb{I}_{\{\#b_*>k\}}+\frac{\theta+\alpha\#\pi'}{(n-1)k+\theta}\mathbb{I}_{\{\#b_*=k\}}\right]\\
&=&\varepsilon^{(n-1)}_{\alpha,\theta}(\pi').
\end{eqnarray*}
This completes the proof.
\end{proof}
  
Though we do not pursue it in detail, we conclude this section with the observation that
exchangeable divisible partition structures are in correspondence with Kingman's paintbox measures.
In particular, the $\varepsilon_{\alpha,\theta}$-family of measures in \eqref{eq:even fidi} is in correspondence with the $\text{Poisson--Dirichlet}(\alpha,\theta)$ measures, for all $k=1,2,\ldots$.

\begin{theorem}
Under the above deletion scheme, exchangeable divisible partition structures are in one-to-one correspondence with Kingman's paintbox measures.
\end{theorem}
\begin{proof}
We need only sketch the proof since the arguments follow by Theorems \ref{thm:CRP even} and \ref{prop:even partition structures}
and Kingman's paintbox representation.  Given a collection $(\Pi_n,\,n\geq1)$ of exchangeable $k$-divisible partitions that is
consistent in distribution under the above divisible deletion scheme, we can obtain a collection $(\Pi_n^*,\,n\geq1)$ of
exchangeable partitions of $(\partitionsn,\,n\geq1)$ by ``deflating'' each block by a factor of $k$ and choosing a representative
element $1,\ldots,n$ of each group of size $k$ within each block.  The result will be an exchangeable partition of $[n]$, which
must obey one of Kingman's paintbox distributions.  The rest now follows by analogous argument to previous theorems.  
\end{proof}

\subsection{Divisible Markov structures}\label{section:even Markov structures}

For $\alpha>0$ and $m\in\mathbb{N}$, the $\PY(-\alpha,m\alpha)$ distribution
determines a probability measure on the subspace $\partitionsNm$ of partitions of
$\mathbb{N}$ having at most $m$ blocks.  Previously, Crane \cite{Crane2011a} introduced an exchangeable
family of Markov chains on $\partitionsNm$ with
marginal transition probabilities
\begin{equation}\label{eq:fidi tps}
p^{(n)}_{\alpha,m}(\pi,\pi')=m^{\downarrow\#\pi'}\prod_{b\in\pi}\frac{\prod_{b'\in\pi'}(\alpha/m)^{\uparrow\#(b\cap
b')}}{\alpha^{\uparrow\#b}},\quad\pi,\pi'\in\partitionsnm,\end{equation}
where $m^{\downarrow n}:=m(m-1)\cdots(m-n+1)$. 
More recently, structural properties of exchangeable Feller processes on $\partitionsNm$ have been characterized in full \cite{Crane2014AOP}.
 The transition probabilities in
\eqref{eq:fidi tps} are reversible with respect to $\varrho^{(n)}_{-\alpha,m\alpha}$ for every $n\in\mathbb{N}$.  
We now extend this family to divisible partitions.

For $\alpha>0$ and
$m\in\mathbb{N}$, the distribution \eqref{eq:even fidi} with parameter
$(-\alpha k,m\alpha k)$ is
\begin{equation}\label{eq:even fidi alpha m}
\varepsilon_{-\alpha k,m\alpha
k}^{(n)}(\pi)=m^{\downarrow\#\pi}\frac{n!}{(nk)!}\frac{\prod_{b\in\pi}\alpha^{\uparrow\#b/k}\frac{\#b!}{(\#b/k)!}}{(m\alpha)^{\uparrow
n}}.\end{equation}
Let $p_{\alpha,m}^{(n)}$ denote the transition probabilities in \eqref{eq:fidi tps}, and let $\evenpartitionsnkm$ be the subset of $k$-divisible partitions of
$[nk]$ with at most $m$ blocks.  We describe a Markovian transition procedure on
$\evenpartitionsnkm$ as follows.  Fix $\pi\in\evenpartitionsnkm$.
\begin{itemize}
	\item[(i)] Independently, for each $b\in\pi$, randomly partition $b$ into $\#b/k$ groups of size $k$ according to the uniform distribution on such partitions of $b$. Label the groups uniquely in $[n]$ to obtain a collection  of groups $\{g_1,\ldots,g_n\}$.
	\item[(ii)] Given $\{g_1,\ldots,g_n\}$ from (i), let $\pi^*$ denote the partition of $\{g_1,\ldots,g_n\}$ obtained by regarding each group as a single element in $\pi^*$, and generate $\Pi''\sim p^{(n)}_{\alpha,m}(\pi^*,\cdot)$, as in \eqref{eq:fidi tps}.
	\item[(iii)] Given $\Pi''=\pi''$, obtain the next state $\pi'\in\evenpartitionsnkm$ by replacing each $g_i$ in $\pi''$ with the group of $k$ elements it represents from (i).
\end{itemize}
\begin{ex}
To illustrate the above transition procedure, let $n=3$, $k=2$, and $\pi=1246/35$.  We generate the transition $\pi\mapsto\pi'$ as follows.
\begin{itemize}
	\item[(i)] We randomly partition the blocks of $\pi$ into sub-blocks of size $k$ and assign labels $1,2,3$, e.g., $g_1=14$, $g_2=26$, and $g_3=35$;
	\item[(ii)] The above procedure yields a partition $\pi^*=12/3$, from which we generate $\Pi''$ according to $p_{\alpha,m}^{(n)}(\pi^*,\cdot)$, say $\pi''=13/2$.
	\item[(iii)] We obtain $\pi'$ by substituting $g_i$ for $i$ in $\pi''$, i.e., $\pi'=g_1g_3/g_2=1345/26$.
\end{itemize}
As in Example \ref{ex:CRP}, there is more than one way to generate the transition $\pi\mapsto\pi'$.  
We derive the transition probability in Theorem \ref{thm:even tps}.

\end{ex}

In the following theorem, we write $\pi\wedge\pi'$ to denote the usual {\em meet} of $\pi$ and $\pi'$, that is,
\[\pi\wedge\pi':=\{B_i\cap B_j':\,B_i\in\pi,\,B'_j\in\pi'\}\setminus\{\emptyset\}.\]
\begin{theorem}\label{thm:even tps} The finite-dimensional transition probabilities of the transition procedure in (i)-(iii) are
\begin{equation}\label{eq:even tps}
\mathcal{E}^{(n)}_{\alpha,m}(\pi,\pi')=m^{\downarrow\#\pi'}\prod_{b\in\pi}\left[\frac{(\#b/k)!}{\#b!}\frac{1}{\alpha^{\uparrow\#b/k}}\prod_{b'\in\pi'}\frac{\#(b\cap
b')!}{\left[\frac{\#(b\cap b')}{k}\right]!}(\alpha/m)^{\uparrow{\#(b\cap
b')/k}}\right],\end{equation}
for $\pi,\pi'\in\evenpartitionsnkm$ satisfying $\pi\wedge\pi'\in\evenpartitionsnk$.
Moreover, for each $n\in\mathbb{N}$, $\mathcal{E}^{(n)}_{\alpha,m}$ is reversible with respect to $\varepsilon^{(n)}_{-\alpha k,m\alpha k}$ in \eqref{eq:even fidi} and is exchangeable with respect to the symmetric group on $[nk]$.
\end{theorem}

\begin{proof}
Fix $n\in\mathbb{N}$ and let $\pi,\pi'\in\evenpartitionsnkm$ satisfy $\pi\wedge\pi'\in\evenpartitionsnk$.  According to the
transition procedure, we first group elements within each $b\in\pi$
together in groups of size $k$.  There are
\[\frac{\#b!}{(k!)^{\#b/k}}\]
ways to group elements and subsequently label them uniquely, each equally likely.  Given a block $b\in\pi$ and the new partition $\pi'$, there are
\[\frac{(\#b/k)!}{(k!)^{\#b/k}}\prod_{b'\in\pi'}\frac{\#(b\cap
b')!}{\left[\frac{\#(b\cap b')}{k}\right]!}\]
sets of labeled groups of $b$ for which a transition $\pi\mapsto\pi'$ is permissible.  (For each $b'\in\pi'$, there are 
\[\frac{\#(b\cap b')!}{(k!)^{\#(b\cap b')/k}}\]
labeled partitions of the elements of $b\cap b'$ into blocks of size $k$.  Dividing this number by $\left[\frac{\#(b\cap b)}{k}\right]!$ gives the number of unlabeled partitions of $b\cap b'$ into $\#(b\cap b')/k$ blocks of size $k$.  Given a partition of $b$, there are $(\#b/k)!$ ways to label the blocks.)
Each permissible partition of $b$ has probability 
\[\frac{(k!)^{\#b/k}}{\#b!}.\]
Multiplication of the number of groupings by their probabilities gives a total
factor of
\begin{equation}\label{eq:factor1}\prod_{b\in\pi}\frac{(\#b/k)!}{\#b!}\prod_{b'\in\pi'}\frac{\#(b\cap
b')!}{\left[\frac{\#(b\cap b')}{k}\right]!}.\end{equation}
Given a partition $\pi^*$ of $n$ groups of size $k$, we choose $\pi^{**}$ from the
transition probabilities of \eqref{eq:fidi tps},
 \begin{equation}\label{eq:factor2}m^{\downarrow\#\pi'}\prod_{b\in\pi}\frac{\prod_{b'\in\pi'}(\alpha/m)^{\uparrow(\#(b\cap
b')/k)}}{\alpha^{\uparrow\#b/k}}.\end{equation}
 Multiplying \eqref{eq:factor1} and \eqref{eq:factor2} gives \eqref{eq:even tps}.
 Reversibility is clear by checking the detailed balanced condition, and exchangeability is clear by inspection.
This completes the proof.

 \end{proof}

\subsection{Divisible Markovian deletion}\label{section:Markovian deletion}
In the following deletion scheme, we define $\sigma(\pi)$ as the image of $\pi\in\partitionsn$ by a permutation $\sigma:[n]\rightarrow[n]$, that is, $i$ and $j$ are in the same block of $\sigma(\pi)$ if and only if $\sigma^{-1}(i)$ and $\sigma^{-1}(j)$ are in the same block of $\pi$.

For $n\in\mathbb{N}$, let $\Pi^{n+1}=(\Pi^{n+1}_1,\Pi^{n+1}_2,\ldots)$ be a Markov chain on $\mathcal{P}^{(m)}_{[(n+1)k]\,\mid\, k}$.  Given $\Pi^{n+1}=(\pi^{n+1}_j,\,j\geq1)$, we obtain the sequence $\Pi^n=(\pi^n_j,\,j\geq1)$ in $\evenpartitionsnkm$ as follows.
 \begin{itemize}
 	\item[(i)] Obtain $\pi^n_1$ from $\pi^{n+1}_1$ by the divisible deletion scheme in Section \ref{section:even partition structures}.  Let $\sigma_1$ be the permutation, called the {\em displacement}, generated in \eqref{eq:even displacement}.
 	\item[(ii)] For $j\geq1$, given $(\Pi^{n+1}_1,\ldots,\Pi^{n+1}_{j+1})=(\pi^{n+1}_1,\ldots,\pi^{n+1}_{j+1})$, $(\Pi^n_1,\ldots,\Pi^n_j)=(\pi^n_1,\ldots,\pi^n_j)$, and displacements $\sigma_1,\ldots,\sigma_j$, we put $\sigma^{(j)}=\sigma_j\circ\cdots\circ\sigma_1$, denote $\pi^*=\sigma^{(j)}(\pi^{n+1}_j)$ and $\pi'=\sigma^{(j)}(\pi^{n+1}_{j+1})$, and let $b_*\in\pi^*$, $b'_*\in\pi'$ be the blocks containing element $nk+1$.
 	\begin{itemize}
 		\item[(a)] Sequentially, for $i=k,\ldots,2$, choose $u^{(i)}$ uniformly from
 		\[ b_*\cap b_*'\cap([nk-k+i]\setminus\{nk-k+1,u^{(k)},\ldots,u^{(i+1)}\})\] and put
 		\begin{equation}\label{eq:even Markov displacement}
 		\sigma_{j+1}=\begin{pmatrix} nk-k+2 \\ u^{(2)}\end{pmatrix} \cdots \begin{pmatrix} nk\\
 	u^{(k)}
 		\end{pmatrix}.\end{equation}
 		\item[(b)] Let $\pi''=\sigma_{j+1}(\pi')$ and obtain $\pi^n_{j+1}$ from $\pi''$ by deleting $\{nk+1,\ldots,nk+k\}$.
 	\end{itemize}
 \end{itemize}
\begin{remark}
The relabeling in (ii) by composing the displacements ensures that the elements $1,2,\ldots$ are consistently labeled in the restricted sequence.  This is needed
because of the first step of the transition procedure, whereby individuals are randomly grouped into sub-blocks of size $k$.
\end{remark}

\begin{remark}
For each $n=1,2,\ldots$, the Markovian deletion operation  $\Pi^{n+1}\longrightarrow\Pi^{n}$ is more than an independent application of the divisible deletion scheme from Section \ref{section:even partition structures} at each time.
For each $t=1,2,\ldots$, $(\Pi^{n+1}_t,\Pi^n_t)$ is marginally distributed as a $k$-divisible partition structure, but Step (ii) of the Markovian deletion scheme incorporates dependence among the deletions across time.
\end{remark}

\begin{definition}
We call a collection $\Pi=(\Pi^1,\Pi^2,\ldots)$ of Markov chains an $(\varepsilon,\mathcal{E})$-{\em reversible Markov structure} if, for each $n\in\mathbb{N}$, $\Pi^n$ is an $(\varepsilon^{(n)},\mathcal{E}^{(n)})$-Markov chain that is reversible with respect to $\varepsilon^{(n)}$, and $\Pi'$ obtained by applying the above deletion scheme to $\Pi^{n+1}$ is an $(\varepsilon^{(n)},\mathcal{E}^{(n)})$-Markov chain.
\end{definition}
 
 \begin{theorem}\label{thm:even Markov structure} For each $n\in\mathbb{N}$, let $\Pi^n$ be an $(\varepsilon^{(n)}_{-\alpha k, m\alpha k},\mathcal{E}^{(n)}_{\alpha,m})$-Markov chain on $\evenpartitionsnkm$.  Then $\Pi=(\Pi^n,\,n\in\mathbb{N})$ is an $(\varepsilon_{-\alpha k,m\alpha k},\mathcal{E}_{\alpha,m})$-reversible Markov structure.
 \end{theorem}
 \begin{proof}
 Fix $n\in\mathbb{N}$ and let $\Pi^{n+1}=(\Pi^{n+1}_j,\,j\geq1)$ be an $(\varepsilon^{(n+1)}_{-\alpha k,m\alpha k},\mathcal{E}^{(n+1)}_{\alpha,m})$-Markov chain.  Let $\Pi^n=(\Pi^n_j,\,j\geq1)$ be obtained from $\Pi^{n+1}$ by the above deletion scheme.  By Theorem \ref{prop:even partition structures} and step (i) of the deletion process, $\Pi^n_1\sim\varepsilon^{(n)}_{-\alpha k,m\alpha k}$.  We induct on $j$ to show that $\Pi^n=(\Pi^n_j,\,j\geq1)$ is an $(\varepsilon^{(n)}_{-\alpha k,m\alpha k},\mathcal{E}^{(n)}_{\alpha,m})$-Markov chain.  
 
 For $j\geq1$, assume the marginal law of $\Pi^n_j$ is $\varepsilon^{(n)}_{-\alpha k,m\alpha k}$ and let $\sigma_1,\ldots,\sigma_j$ be the displacement permutations used in the Markovian deletion scheme up to step $j$.  By Theorem \ref{thm:CRP even}, $\sigma^{(j)}(\Pi^{n+1}_j)\sim\varepsilon^{(n)}_{-\alpha k,m\alpha k}$ and, by Theorem \ref{thm:even tps}, the conditional distribution of $\sigma^{(j)}(\Pi^{n+1}_{j+1})$ given $\sigma^{(j)}(\Pi^{n+1}_j)$ is $\mathcal{E}^{(n+1)}_{\alpha,m}(\sigma^{(j)}(\Pi^{n+1}_j),\cdot)$.  Given $(\sigma^{(j)}(\Pi^{n+1}_j),\sigma^{(j)}(\Pi^{n+1}_{j+1}))=(\pi,\pi')$ and $\Pi^n_j=\pi^*$, let $\pi''\in\evenpartitionsnkm$ be such that $\mathcal{E}^{(n)}_{\alpha,m}(\pi^*,\pi'')>0$ and there is a positive probability of obtaining $\pi''$ from $\pi'$ in the divisible Markovian deletion procedure.  Let $b_*\in\pi$, $b_*'\in\pi'$ be the blocks containing $nk+1$. We have

\begin{eqnarray}
\lefteqn{\frac{\mathcal{E}^{(n+1)}_{\alpha,m}(\pi,\pi')}{\mathcal{E}^{(n)}_{\alpha,m}(\pi^*,\pi'')}=}\notag\\
&&\frac{\#b_*/k}{(\#b_*)^{\downarrow k}}\frac{1}{\alpha+\#b_*/k-1}\left[\left(\frac{\alpha}{m}+\frac{\#(b_*\cap b'_*)}{k}-1\right)\frac{(\#(b_*\cap b'_*))^{\downarrow k}}{\left[\frac{\#(b_*\cap b'_*)}{k}\right]}\mathbb{I}_{\{\#b'_*>k\}}+k!\frac{(m-\#\pi')\alpha}{m}\mathbb{I}_{\{\#b'_*=k\}}\right].\label{eq:conditional tp even}\end{eqnarray}
In step (ii)(a), $\sigma_{j+1}$ in \eqref{eq:even Markov displacement} has probability $(\#(b_*\cap b'_*))/((\#(b_*\cap b'_*))^{\downarrow k})$.

Now, for every $\pi'$ with $\mathcal{E}^{(n)}_{\alpha,m}(\pi,\pi')>0$, the block sizes of $\pi'$ are the same as $\pi''$ except for the block $b'_*\in\pi'$, which contains $k$ more elements than its corresponding block in $\pi''$.  Given $b'_*$, there are $(\#b_*)^{\downarrow k}/\#b_*$ partitions $\pi''$ from which $\pi'$ can be obtained by the even Markovian deletion process.  Each $\pi''$ corresponds to a unique displacement $\sigma_{j+1}$ in \eqref{eq:even Markov displacement}.  By the law of cases, we have
\begin{eqnarray*}
\lefteqn{\mathbb{P}\{\Pi^{n}_{j+1}=\pi''\mid \Pi^n_{j}=\pi^*,\,\Pi^{n+1}_j=\pi\}=}\\
&=&\sum_{\pi'}\mathbb{P}\{\Pi^{n}_{j+1}=\pi''\mid \Pi^n_{j}=\pi^*,\,(\Pi^{n+1}_j,\Pi^{n+1}_{j+1})=(\pi,\pi')\}\,\mathcal{E}^{(n+1)}_{\alpha,m}(\pi,\pi')\\
&=&{\mathcal{E}^{(n)}_{\alpha,m}(\pi^*,\pi'')}\sum_{b'_*}\sum_{\sigma_{j+1}}\frac{\#b_*}{(\#b_*)^{\downarrow k}}\frac{1}{\alpha+\#b_*/k-1}\left[\left(\frac{\alpha}{m}+\frac{\#(b_*\cap b'_*)}{k}-1\right)\mathbb{I}_{\{\#b'_*>k\}}+\frac{(m-\#\pi')\alpha}{m}\mathbb{I}_{\{\#b'_*=k\}}\right]\\
&=&\mathcal{E}^{(n)}_{\alpha,m}(\pi^*,\pi'')\sum_{b'_*}\frac{1}{\alpha+\#b_*/k-1}\left[\left(\frac{\alpha}{m}+\frac{\#(b_*\cap b'_*)}{k}-1\right)\mathbb{I}_{\{\#b'_*>k\}}+\frac{(m-\#\pi')\alpha}{m}\mathbb{I}_{\{\#b'_*=k\}}\right]\\
&=&\mathcal{E}^{(n)}_{\alpha,m}(\pi^*,\pi'').
\end{eqnarray*}
By the induction hypothesis and reversibility (Theorem \ref{thm:even tps}), the unconditional law of $\Pi^n_{j+1}$ is $\varepsilon^{(n)}_{-\alpha k,m\alpha k}$ and $\Pi^n$ is an $(\varepsilon^{(n)}_{-\alpha k,m\alpha k},\mathcal{E}^{(n)}_{\alpha,m})$-Markov chain.
This completes the proof.
\end{proof}

\section{Concluding remarks}\label{section:remarks}
\subsection{Divisible paintbox partitions}
We have treated the special case of Kingman's paintbox process directed by
$\nu=\text{Poisson--Dirichlet}(\alpha,\theta)$; however, we can
apply the above descriptions more generally to generate both divisible partition structures
and divisible Markov structures associated to any probability measure on the ranked-simplex.
In this case, we replace the Ewens--Pitman measure by a paintbox measure $\varrho_{\nu}$
in the above statements, but reversibility of the resulting family of
Markov structures is not guaranteed.  The Markov partition structures for the $k=1$ case have been
studied in \cite{Crane2011a}; the $k>1$ case follows by modifying the work in \cite{Crane2011a} according to the program in Section \ref{section:even Markov structures}.

\subsection{Divisible random permutations}
There is a well-known connection between exchangeable random partitions and random permutations whose
distribution depends only on their cycle sizes.  The combinatorial identity \eqref{eq:identity} relates 
divisible permutations and divisible partitions in the usual way: each partition of $[nk]$ with blocks of size
$(n_1,\ldots,n_m)$ corresponds to $\prod_{j=1}^m(n_j-1)!$ permutations of $[nk]$ with the corresponding
cycle sizes.  The usual approach of sampling uniformly from the subset of permutations corresponding to a given
random partition yields a distribution on what we call {\em $k$-divisible permutations}. 
Divisible permutations appear in combinatorics, e.g., Wilf \cite{WilfGF}, and the generating function \eqref{eq:identity} for $k$-divisible 
permutations generates known integer sequences, e.g., \cite{OEIS}:A001818 and \cite{OEIS}:A178575.

\subsection{Neutral partition structures}

We can also consider partitions of a marked population
 $[n]^*$ consisting of $nk$ individuals, with $n$ individuals of each type $j=1,\ldots,k$.  We write $1^{(j)},\ldots,n^{(j)}$ to denote the individuals
of type $j=1,\ldots,k$.  We call a subset of $[n]^*$ {\em neutral} if it contains
an equal number of elements of each type.  A partition of $[n]^*$ is neutral if each of
its blocks is.  Necessarily, any neutral partition is $k$-divisible.

As in the $k$-divisible case, $k=1$ corresponds to ordinary set partitions; however,
by adding more structure to the partitions, the combinatorial arguments for neutral
partitions are more straightforward.  The only difference from the $k$-divisible case
is that our seating rule is specified to preserve neutrality, not just divisibility.  For example,
the Chinese restaurant process in the neutral case corresponds to the following
seating rule.

\vspace{2mm}

\noindent{\bf Neutral random seating rule}
\begin{itemize}
        \item[(1)] The first $k$ individuals are seated at the same table,
$\Pi_1:=\{1^{(1)},\ldots,1^{(k)}\}$.
        \item[(2)] After $nk$ arrivals are seated according to $\Pi_n$, the next $k$
individuals $(n+1)^{(1)},\ldots,(n+1)^{(k)}$ are seated randomly as follows:
        \begin{itemize}
                \item[(a.)] independently for each $i=2,\ldots,k$, $(n+1)^{(i)}$ chooses $u^{(i)}$
uniformly among $1^{(i)},\ldots,(n+1)^{(i)}$ and displaces $u^{(i)}$ in $\Pi_n$;
                \item[(b.)] the displaced group $((n+1)^{(1)},u^{(2)},\ldots,u^{(k)})$ is treated
as a single unit and randomly sits at table $b\in\Pi_n\cup\{\emptyset\}$ according to $\crp(nk,\alpha,\theta)$.
        \end{itemize}
\end{itemize}
The finite-dimensional distributions on neutral partitions of $[n]^*$ generated by the above seating rule are
   \begin{equation}\label{eq:balanced fidi}
  \nu^{(n)}_{\alpha,\theta}(\pi)=\frac{(\theta/\alpha)^{\uparrow\#\pi}}{(\theta/k)^{\uparrow
n}}\frac{\prod_{b\in\pi}-(-\alpha/k)^{\uparrow\#b/k}[(\#b/k)!]^{k-1}}{(n!)^{k-1}}.
   \end{equation}  
The distribution in \eqref{eq:balanced fidi} is exchangeable with respect to permutations of $[n]^*$ that permute
only elements with the same type.  We obtain the deletion rule by reversing the Chinese restaurant seating
rule, as we have for the $k$-divisible deletion rule. 

Neutral Markov structures are collections of pairs $\{(\nu^{(n)},\mathcal{N}^{(n)})\}_{n\geq1}$ so that
$\{\nu^{(n)}\}_{n\geq1}$ is an exchangeable neutral partition structure and $\{\mathcal{N}^{(n)}\}_{n\geq1}$
is a collection of exchangeable Markovian transition probabilities consistent under an operation of neutral
Markovian deletion, which is analogous to the divisible Markovian deletion scheme.  In the neutral setting,
we obtain finite-dimensional transition probabilities corresponding to the Ewens--Pitman Markov chain:
\begin{equation}\label{eq:neutral tps}
\mathcal{N}^{(n)}_{\alpha,m}(\pi,\pi')=m^{\downarrow\#\pi'}\prod_{b\in\pi}\left[\frac{\prod_{b'\in\pi'}(\alpha/m)^{\uparrow(\#(b\cap b')/k)}[(\#(b\cap b')/k)!]^{k-1}}{[(\#b/k)!]^{k-1}\alpha^{\uparrow\#b/k}}\right],\end{equation}
provided $\pi\wedge\pi'$ is a neutral partition.  The transition probability in \eqref{eq:neutral tps} is reversible with respect to $\nu^{(n)}_{-\alpha k,m\alpha k}(\cdot)$ and is exchangeable with respect to permutations of $[n]^*$ that preserve neutrality.  Connections between exchangeable
neutral partition structures and Kingman's paintbox measures are analogous to those for divisible partition
structures and follow by similar arguments.

\acks
Harry Crane is partially supported by U.S.\ NSF DMS-1308899 and U.S.\ NSA H98230-13-1-0299 grants.


\begin{thebibliography}{10}

\bibitem{BaileyAssociationSchemes}
R.~A. Bailey.
\newblock {\em Association Schemes: Designed Experiments, Algebra and
  Combinatorics}, volume~84 of {\em Cambridge Studies in Advanced Mathematics}.
\newblock Cambridge University Press, Cambridge, 2004.

\bibitem{Bertoin2006}
J.~Bertoin.
\newblock {\em Random fragmentation and coagulation processes}, volume 102 of
  {\em Cambridge Studies in Advanced Mathematics}.
\newblock Cambridge University Press, Cambridge, 2006.

\bibitem{Bertoin2008}
J.~Bertoin.
\newblock Two-parameter {P}oisson-{D}irichlet measures and reversible
  exchangeable fragmentation-coalescence processes.
\newblock {\em Combin. Probab. Comput.}, 17(3):329--337, 2008.

\bibitem{Crane2011a}
H.~Crane.
\newblock A consistent {M}arkov partition process generated from the paintbox
  process.
\newblock {\em J. Appl. Probab.}, 43(3):778--791, 2011.

\bibitem{Crane2014cluster}
H.~Crane.
\newblock Clustering from categorical data sequences.
\newblock {\em Journal of the American Statistical Association}, in press, 2014.

\bibitem{Crane2013b}
H.~Crane.
\newblock Permanental {P}artition {M}odels and {M}arkovian {G}ibbs
  {S}tructures.
\newblock {\em Journal of Statistical Physics}, 153(4):698--726, 2013.

\bibitem{Crane2013LAA}
H.~Crane.
\newblock Some algebraic identities for $\alpha$-permanent.
\newblock {\em Linear Algebra and Its Applications}, 439(11):3445-3459, 2013.

\bibitem{Crane2014AOP}
H.~Crane.
\newblock The cut-and-paste process.
\newblock {\em Annals of Probability}, 42(5):1952--1979, 2014.

\bibitem{DonnellyGrimmett1993}
P.~Donnelly and G.~Grimmett.
\newblock On the asymptotic distribution of large prime factors.
\newblock {\em J. London Math. Soc}, 47:395--404, 1993.

\bibitem{EfronThisted1976}
B.~Efron and R.~Thisted.
\newblock Estimating the number of unseen species: How many words did
  shakespeare know?
\newblock {\em Biometrika}, 63:435--447, 1976.

\bibitem{EfronThisted1987}
B.~Efron and R.~Thisted.
\newblock Did shakespeare write a newly discovered poem?
\newblock {\em Biometrika}, 74:445--455, 1987.

\bibitem{Ewens1972}
W.~J. Ewens.
\newblock The sampling theory of selectively neutral alleles.
\newblock {\em Theoret. Population Biology}, 3:87--112, 1972.

\bibitem{FisherCorbetWilliams1943}
R.~Fisher, A.~Corbet, and C.~Williams.
\newblock The relation between the number of species and the number of
  individuals in a random sample of an animal population.
\newblock {\em J. Animal Ecology}, 12:42--58, 1943.

\bibitem{GnedinHaulkPitman}
A.~Gnedin, C.~Haulk, and J.~Pitman.
\newblock Characterizations of exchangeable partitions and random discrete
  distributions by deletion properties.
\newblock In {\em Probability and mathematical genetics}, volume 378 of {\em
  London Math. Soc. Lecture Note Ser.}, pages 264--298. Cambridge Univ. Press,
  Cambridge, 2010.

\bibitem{Kingman1978a}
J.~F.~C. Kingman.
\newblock Random partitions in population genetics.
\newblock {\em Proc. Roy. Soc. London Ser. A}, 361(1704):1--20, 1978.

\bibitem{Kingman1978b}
J.~F.~C. Kingman.
\newblock The representation of partition structures.
\newblock {\em J. London Math. Soc. (2)}, 18(2):374--380, 1978.


\bibitem{McCullaghYang2008}
P.~McCullagh and J.~Yang.
\newblock How many clusters?
\newblock {\em Bayesian Anal.}, 3(1):101--120, 2008.

\bibitem{Pitman2005}
J.~Pitman.
\newblock {\em Combinatorial stochastic processes}, volume 1875 of {\em Lecture
  Notes in Mathematics}.
\newblock Springer-Verlag, Berlin, 2006.
\newblock Lectures from the 32nd Summer School on Probability Theory held in
  Saint-Flour, July 7--24, 2002, With a foreword by Jean Picard.

\bibitem{PitmanYor1997}
J.~Pitman and M.~Yor.
\newblock The two-parameter {P}oisson-{D}irichlet distribution derived from a
  stable subordinator.
\newblock {\em Ann. Probab.}, 25(2):855--900, 1997.

\bibitem{OEIS}
N.~Sloane.
\newblock Online {E}ncyclopedia of {I}nteger {S}equences.
\newblock {\em Published electronically at http://www.oeis.org/}.

\bibitem{WilfGF}
H.~S. Wilf.
\newblock {\em {g}eneratingfunctionology, Third Edition}.
\newblock AK Peters/CRC Press, 2005.

\end{thebibliography}
\end{document}